    \documentclass[11pt]{article}
\newcommand{\ddate}{13 August 2013}

\date{\ddate}

\usepackage{amssymb}

\usepackage[matrix,arrow,curve]{xy}

\usepackage{epsfig}

\usepackage{color}

\usepackage{psfrag}

\usepackage{amsmath}

\numberwithin{equation}{section}

\numberwithin{figure}{section}

\usepackage{theorem}

\newtheorem{thm}{Theorem}[section]

\newtheorem{Theorem}[thm]{Theorem}

\newtheorem{Lemma}[thm]{Lemma}

\newtheorem{Proposition}[thm]{Proposition}

\newtheorem{Corollary}[thm]{Corollary}

{\theorembodyfont{\rmfamily} \theoremstyle{plain}

\newtheorem{rem}[thm]{Remark}

\newtheorem{Remark}[thm]{Remark}

\newtheorem{Example}[thm]{Example}

\newtheorem{definition}[thm]{Definition}

\newcommand{\preu}{\noindent{\sc Proof: \ }}

\newcommand{\fl}[1]{\buildrel{#1}\over{\longrightarrow}}

\newcommand{\cqfd}{\unskip\kern 6pt\penalty 500
\raise -2pt\hbox{\vrule\vbox to10pt{\hrule width
 8pt\vfill\hrule}\vrule}\smallskip}

\newcommand{\bbd}{{\mathbb{D}}}

\newcommand{\bbs}{{\mathbb{S}}}

\newcommand{\bbp}{{\mathbb{P}}}

\newcommand{\bbe}{{\mathbb{E}}}

\newcommand{\bbr}{{\mathbb{R}}}

\newcommand{\bbc}{{\mathbb{C}}}

\newcommand{\bbz}{{\mathbb{Z}}}

\newcommand{\bbl}{{\mathbb{L}}}

\newcommand{\calc}{{\mathcal C}}

\newcommand{\cald}{{\mathcal D}}

\newcommand{\cale}{{\mathcal E}}

\newcommand{\calg}{{\mathcal G}}

\newcommand{\calh}{{\mathcal H}}

\newcommand{\calj}{{\mathcal J}}

\newcommand{\calm}{{\mathcal M}}

\newcommand{\caln}{{\mathcal N}}

\newcommand{\calp}{{\mathcal P}}

\newcommand{\cals}{{\mathcal S}}

\newcommand{\pcirc}{\kern .7pt {\scriptstyle \circ} \kern 1pt}

\newcounter{exo}

\newcommand{\mun}{{-1}}

\newcommand{\sk}[1]{\vskip #1 mm}

\newcommand{\llangle}[2]{\langle #1 ,#2 \rangle}

\newcommand{\onto}{\to\kern-6.5pt\to}

\newcommand{\algt}{\mathfrak{t}}

\newcommand{\proref}[1]{Proposition~\ref{#1}}

\newcommand{\remref}[1]{Remark~\ref{#1}}

\newcommand{\lemref}[1]{Lemma~\ref{#1}}

\newcommand{\corref}[1]{Corollary~\ref{#1}}

\newcommand{\thref}[1]{Theorem~\ref{#1}}

\newcommand{\exref}[1]{Example~\ref{#1}}

\newcommand{\simple}{simple Hamiltonian}

\renewcommand{\:}{\colon}

\newcommand{\grad}{{\rm Grad\,}}

\newcommand{\nua}[2]{\caln^{#1}_{#2}}

\newcommand{\pol}[1]{{\caln}_{#1}}

\newcommand{\rmcr}{\ensuremath{\bbr_{\scriptscriptstyle\nearrow}^n}}

\title{Simple Hamiltonian manifolds}

\author{Jean-Claude HAUSMANN and Tara HOLM}

\date{\ddate}
\begin{document}%aaa

\maketitle

\begin{abstract}
A {\bf\simple\ manifold} is a compact connected symplectic manifold equipped with a Hamiltonian action of a 
torus $T$ with moment map $\Phi: M\to \algt^*$,  such that $M^T$ has exactly two connected components, 
denoted $M_0$ and $M_1$.  We study the differential and symplectic geometry of \simple\ manifolds, 
including a large number of examples.
\end{abstract}

\tableofcontents

\section{Introduction}\label{hamedge}
%Definitions and examples

Let $M$ be a compact connected symplectic manifold equipped with a Hamiltonian action of a torus $T=(S^1)^n$, and let 
$\Phi:M\to\algt^*$ denote the moment map.   The celebrated Atiyah Guillemin-Sternberg convexity theorem states 
that the image of the moment map $\Phi$ is the convex hull of the image of the fixed points, $\Phi(M^T)$.
This polytope is a single point, that is, the moment map is constant, if and only if the action is trivial.  So long as
the action is non-trivial, this polytope $\Phi(M)$ must have at least two extreme points.  In this paper, we consider
the simplest non-trivial case, when $M^T$ has exactly two components, and so $\Phi(M)$ is a $1$-dimensional polytope.

\begin{definition}\label{def:simple}
A {\bf\simple\ manifold} is a compact connected symplectic manifold equipped with a Hamiltonian action of a 
torus $T$ with moment map $\Phi: M\to \algt^*$,  such that $M^T$ has exactly two connected components, 
denoted $M_0$ and $M_1$.  
\end{definition}

As noted above, a simple manifold has
the minimum possible number of fixed components. 
We describe a simple Hamiltonian manifold by the triple $(M,M_0,M_1)$, and let $2m_i$ and $2m$ be the dimensions of 
$M_i$ and $M$ respectively, and set $2r_i={\rm codim\,}M_i=2m-2m_i$.  As a consequence of some basic
results in equivariant symplectic geometry, the torus action on a simple
manifold necessarily factors into a trivial action and a {\bf residual} effective circle action (Lemma~\ref{character}).
Thus, our results hold for torus actions, but generally require verification only for the residual circle action.

In what follows, we explore the geometry associated to \simple\ manifolds.  We establish the basic topology
of a \simple\ manifold, using the moment map as the key tool, in Section~\ref{se:prelim}.  This is where we discuss
the residual circle action (Lemma~\ref{character}).  Then we turn to cohomology
constraints on \simple\ manifolds in Section~\ref{se:coh}.  The residual moment map is a Morse--Bott function
on $M$, and so the cohomology of $M$ is determined from $M_0$ and $M_1$ (Proposition~\ref{HKexseq} and its
Corollaries).  This allows us to deduce relations among $m$, $m_0$, $m_1$, $r_0$ and $r_1$. This section also 
includes comments about how our work relates to several recent papers on this topic.  

In Section~\ref{S.Diffeo}, we study bundles over the $M_i$ and the gauge groups of
these bundles, and prove our first main theorem giving necessary conditions for two \simple\ manifolds to be
$T$-equivariantly diffeomorphic, Theorem~\ref{TDIFF}.  
Next, in Section~\ref{se:case}, we turn to the special case when $M_1$ has 
codimension $2$ in $M$, and characterize $M$ in terms of $M_0$ (Theorem~\ref{T.r1=1}).  
In this special case, we must have that $M_1$ is diffeomorphic to $M_1$ (Corollary~\ref{C-cod2-2}).
In Section~\ref{se:symplectomorphism}, we turn to the classification $M$
up to $T$-equivariant symplectomorphism, with a complete answer in the same special case $r_1=1$ 
(Theorems~\ref{T.symcod1} and \ref{T.S11conn}).  In particular, when $r_0=r_1=1$, then $M_0$
and $M_1$ must be $T$-equivariantly symplectomorphic (Corollary~\ref{T.S1=1}).  Finally, the last 
section of the paper is devoted to examples of polygon spaces.

There has been a flurry of recent work on Hamiltonian $S^1$-manifolds that are in some sense minimal.  Tolman introduces 
Betti number constraints in \cite{T}, and shows that only a finite number of cohomology rings can occur. 
These constraints are explored further in \cite{HT} when the fixed set has exactly two  components, that is the manifold 
is a \simple\ manifold.  The differential geometry of  \simple\ manifolds with minimal Betti numbers is discussed in
\cite{LOS}; this work may be related to our results in Section~\ref{S.Diffeo}.  Another natural hypothesis is that the
circle action be semi-free, as is the case for weight \simple\ manifolds discussed below in Section~\ref{se:prelim}.  
The implications of this hypothesis are developed further in \cite{Gz, TW}.

We now conclude this Introduction with a handful of examples of \simple\ manifolds.

\begin{Example}\label{edgeexfirst}Let $M=\bbc P^n$ with a circle action given by $$
g\cdot [z_0:\dots :z_n]=[gz_0:\dots :gz_k:z_{k+1}:\dots  :z_n] \, ,
$$
for $g\in S^1$.
This is a simple Hamiltonian manifold
$(\bbc P^n,\bbc P^k,\bbc P^{n-k-1})$.
\end{Example}

\begin{Example}
A \simple\ manifold $M$ with $M^T$ discrete
is diffeomorphic to $S^2$. In this case, the moment map is a
Morse function with exactly two critical points, which implies that
$M$ is homeomorphic to a sphere $S^n$. As $M$ is symplectic, it must be
diffeomorphic to $S^2$.
\end{Example}

\begin{Example}\label{eg:sympcut} \em The symplectic cut of a weight bundle. \rm
We may use Lerman's symplectic cuts \cite{Le} to produce a \simple\ manifold from a symplectic 
manifold equipped with a complex vector bundle.
Let $M_0$ be a compact symplectic manifold and let
$\nu_0\:V\to M_0$ be a complex vector bundle of rank $k$. Viewing $S^1\subset \bbc$ as the
unit complex numbers, there is a natural $S^1$-action on this bundle, namely fiberwise complex
multiplication.
We assume that the total space $V$ is equipped with a symplectic form so that
this $S^1$-action is Hamiltonian.
The moment map $\phi\:V\to \bbr$ has only $0$ as a critical value.
Let $M$ be the symplectic cut of $V$ at a regular value $\ell>0$
of $\phi$. This gives a simple Hamiltonian manifold $(M,M_0,M_1)$
with $M_1$ the symplectic reduction of $V$ at $\ell$.
The bundle projection descends to a map $M\to M_0$
with fibre $\bbc P^k$. Thus, $M=\hat\bbp(\nu_0)$, the total space of the
$\bbc P^k$-bundle associated to $\nu_0$ and 
$M_1=\bbp(\nu_0)$, the total space of the
$\bbc P^{k-1}$-bundle associated to $\nu_0$.
The case $k=1$ is described in \cite[Example 5.10]{MDS}.
\end{Example}

\begin{Example}
Let $M=G_k(\bbc^{r})$ be the Grassmannian manifold of complex $k$-planes in $\bbc^r$, endowed with a 
$U(r)$-invariant symplectic form. As a homogeneous space, $M\cong U(r)/(U(k)\times U(r-k))$.  We may endow $M$ with 
a symplectic form by identifying it with the $U(r)$ coadjoint orbit of Hermitian $r\times r$ matrices with eigenvalues
consisting of $k$ ones and $(r-k)$ zeros.  The maximal torus $T$ of diagonal matrices in $U(r)$ acts in a Hamiltonian fashion
on $M$, and we consider the last coordinate circle  of this torus.   Under the identifications we have made, this action 
has moment map
$$
\Phi(A)=a_{r,r}, 
$$
where $A$ is a symmetric matrix and $a_{r,r}$ its bottom right entry.  Then $M$ is a simple Hamiltonian manifold
with moment map image the interval $[0,1]$.
We identify 
$$
M_0 = \left\{ \ \left(\begin{array}{cccc}
 & & & 0\\
 & \mbox{\em \LARGE{B}} & & \vdots \\
 & & & 0 \\
 0 & \cdots & 0 & 0 
 \end{array}\right)\ \right\},
$$ 
where $B$ is a symmetric $(r-1)\times (r-1)$ matrix with eigenvalues
consisting of $k$ ones and $(r-k-1)$ zeros.  Thus,  $M_0\cong G_{k-1}(\bbc^{r-1})$.
The second fixed component is
$$
M_1= \left\{ \ \left(\begin{array}{cccc}
 & & & 0\\
 & \mbox{\em \LARGE{B}} & & \vdots \\
 & & & 0 \\
 0 & \cdots & 0 & 1
 \end{array}\right)\ \right\},
$$
where $B$ is a symmetric $(r-1)\times (r-1)$ matrix with eigenvalues
consisting of $(k-1)$ ones and $(r-k)$ zeros; so $M_1\cong G_k(\bbc^{r-1})$.
The real locus (for complex conjugation) of this simple Hamiltonian manifold is discussed in \cite[Example~5]{HR}.
\end{Example}

\begin{Example}
If $M$ is a \simple\ manifold and $N$ is a
connected compact symplectic manifold, then $M\times N$ is a
\simple\ manifold, where
$g\cdot(x,y)=(gx,y)$ for $g\in T$ and $(x,y)\in M\times N$.
\end{Example}

\begin{Example}\label{EG:SO(5)}
Grassmannian manifold $\tilde G_2(\bbr^{m+2})$ of oriented $2$-planes in $\bbr^{m+2}$. 
See ~\ref{fig:B23} and its legend, describing moment polytopes for $\tilde G_2(\bbr^{5})$ and $\tilde G_2(\bbr^{7})$.
These simple manifolds play an important role in \cite{T,HT,LOS}.

\begin{figure}[ht]
\centering{
\psfrag{A}{(a)}
\psfrag{B}{(b)}
\includegraphics[width=4in]{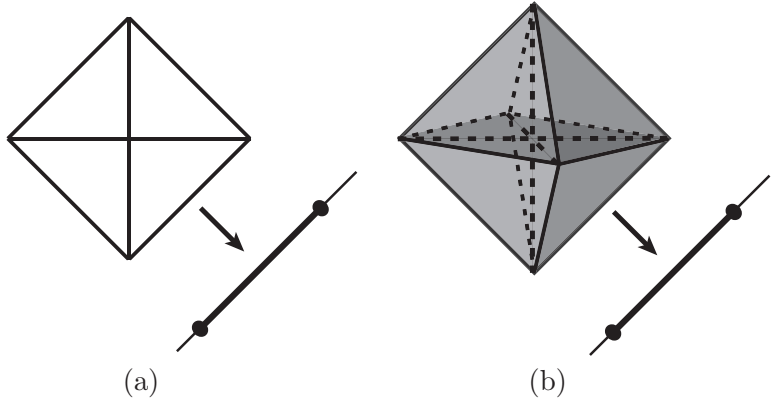}
}
\caption{In Figure (a) there is the $T^2$-moment polytope for $\tilde G_2(\bbr^{5})$, together with a projection 
for an $S^1$-moment map for which this manifold is a \simple\ manifold.  
Both $M_0$ and $M_1$ in this case are diffeomorphic to $\bbp^1$. 
Figure (b) shows the $T^3$-moment polytope for $\tilde G_2(\bbr^{7})$ and a projection for an $S^1$-moment map 
for which this manifold is a \simple\ manifold.  
Both $M_0$ and $M_1$ in this case are diffeomorphic to $\bbp^2$.}
\label{fig:B23}
\end{figure}

\end{Example}

\paragraph{Acknowledgements:} The first author had several  fruitful conversations with
Y.~Karshon and S.~Tolman at the conference on moment maps organized in the Bernoulli 
Center (EPFL) in August 2008.   The second author 
was supported in part by NSF Grant DMS--0835507.  She would like to thank 
Y.~Karshon, D.~McDuff, and S.~Tolman for useful conversations.
Both authors would like to thank the anonymous referees for helpful comments.

%\vfill\pagebreak

\section{Preliminaries}
\label{se:prelim}

Standard properties of moment maps, which may be found in  \cite{Au},
immediately imply the following

\begin{Lemma}\label{edgeeasy}
Let $(M,M_0,M_1)$ be a \simple\ manifold
with moment map $\Phi:M\to\algt^*$. Then
\begin{enumerate}
\renewcommand{\labelenumi}{(\roman{enumi})}
\item the moment polytope $\Delta=\Phi(M)$ is a line segment.
\item $\Phi$ is a Morse-Bott function onto $\Delta$ with exactly two critical values,
namely the endpoints of $\Delta$.
\cqfd
\end{enumerate}
\end{Lemma}

We think of the circle $S^1$ as the complex numbers of norm $1$. The Lie algebra ${\rm Lie}(S^1)$
may then be identified as $i\,\bbr$, with basis vector $2\pi i$. We may use the dual basis to identify
${\rm Lie}(S^1)^*$ with $\bbr$. The group of characters of $T$ is $\hat T=\mathrm{Hom}(T,S^1)$, the set
of smooth homomorphisms. This is isomorphic to the linear maps from $\bbr={\rm Lie}(S^1)^*$
to $\algt^*$ that send  $\bbz$ to the weight lattice. Taking the image of $1$ identifies $\hat T$
with the weight lattice inside $\algt^*$.

\begin{Lemma}\label{character}
Let $(M,M_0,M_1)$ be a \simple\ $T$-manifold
with moment map $\Phi:M\to\algt^*$. Then
there is a unique character $\chi\in\hat T$ such that
\begin{enumerate}
\renewcommand{\labelenumi}{(\roman{enumi})}
\item the $T$-action $\alpha\:T\to {\rm Diff\,}(M)$ is of the form
$\alpha=\bar\alpha\pcirc\chi$ where $\bar\alpha\:S^1\to {\rm Diff\,}(M)$
is an effective action making
$(M,M_0,M_1)$ a simple Hamiltonian $S^1$-manifold.  We call $\bar\alpha$ the {\bf residual action}.
\item The residual action $\bar\alpha$ admits a moment map $\bar\Phi\:M\to\bbr$
such that $\bar\Phi(M_0)=0$ and $\bar\Phi(M_1)>0$.
\end{enumerate}
Moreover, the above character $\chi$, seen as an element of the weight lattice,
is a positive multiple of \, $\Phi(M_1)-\Phi(M_0)$.
\end{Lemma}

\begin{Remark}
The character $\chi$ of Part (ii) of \lemref{character} is the
\textbf{associated character} to the \simple\ manifold $(M,M_0,M_1)$.
The  moment map $\bar\Phi\:M\to\bbr$ is called the \textbf{residual moment map}.
This lemma reduces the classification of \simple\ $T$-manifolds to the case of \simple\ $S^1$-manifolds, for 
effective circle actions. 
\end{Remark}

\preu
As the moment polytope is $1$-dimensional,
$\bar T=T/\ker\alpha$ is a $1$-dimensional torus (see e.g.~\cite[\S~III.2.b]{Au}). 
Choosing an 
identification of $\bar T$ with $S^1$ gives a character $\chi$ 
and a residual action $\bar\alpha$ with moment map $\bar\Phi\:M\to\bbr$
(with $\bar\Phi(M_0)=0$). We denote by $2m_i$ and $2m$ the dimensions of $M_i$ and $M$
and we set $2r_i={\rm codim\,}M_i$.
As $\alpha=\bar\alpha\pcirc\chi$, this implies that $\chi\pcirc\bar\Phi$ is a moment
map for $\alpha$, proving the last statement. 
The uniqueness statement (ii) follows from the fact that the two identifications of $\bar T$
with $S^1$ differ by the sign of $\bar\Phi$.
\cqfd

Let $(M,M_0,M_1)$ be a \simple\ $T$-manifold. 
Recall that $M$ always admits a $T$-invariant almost complex
structure $J$ that is \textbf{$\omega$-compatible}: $J$ is an isometry for $\omega$ and
$\omega(v,Jv)>0$ for all non-zero tangent vectors $v$ to $M$ 
(see, for example, \cite[\S 2.5]{MDS} or \cite[Part~V]{CS}).
Then $\llangle{v}{w}=\omega(v,Jw)$ defines a Riemannian
metric on $M$ and $\llangle{\,}{}+i\omega(,)$ is a $T$-invariant Hermitian metric. 
The space of $T$-invariant $\omega$-compatible almost complex structures on $M$ is denoted by 
$\calj(M,\omega)$ and is contractible (see \cite[Proposition~4.1 and 2.49]{MDS}
or \cite[Prop.~13.1]{CS}). Therefore, choosing $J\in \calj(M,\omega)$
endows the tangent bundle $TM$ with a $U(r)$-structure
whose isomorphism class is well-defined. As the $M_i$ are symplectic submanifolds, the 
normal bundles $\nu_i=TM|_{M_i}/TM_i$ are also Hermitian bundles, with structure group
$U( r_i)$, and these structures are well-defined
up to isomorphism. Observe that $\nu_i$ is isomorphic to the orthogonal complement to $TM_i$
in $TM|_{M_i}$ with respect to the Riemannian metric associated to $J$.

The $U(r_i)$ structure on $\nu_i$ is $T$-invariant, so the bundle $\nu_i$  decomposes
into a Whitney sum of $T$-weight bundles. It follows from Lemmas~\ref{edgeeasy} and~\ref{character}
that the weights which occur are multiples of $\chi$.

\begin{definition}\label{def:weight-bdle}
If $\nu_0$ (or, equivalently,  $\nu_1$) is itself a weight bundle,
we call $M$ a {\bf weight \simple\ manifold.} 
\end{definition}

For instance,
$M$ is a weight \simple\ manifold when ${\rm codim\,}M_0=2$ or
${\rm codim\,}M_1=2$.
 The Grassmannian manifold $\tilde G_2(\bbr^{m+2})$
of Example 1.1.6 is not a weight simple manifold.
Observe that $M$ is a weight \simple\ manifold if and only if the residual action is 
semi-free. By \cite[Proposition~8.1]{HT}, a \simple\ manifold $(M,M_0,M_1)$ with $m=m_0+m_1+1$
is a weight \simple\ manifold unless $\dim M_0=\dim M_1$.

\begin{Remark}\label{R.Char}
In the above discussion the Hermitian bundle $\nu_i$ is the underlying bundle of
a Hermitian bundle $\hat\nu_i$ endowed with a $T$-action. We do not distinguish these two notions
because in the case of interest for us, where $(M,M_0,M_1)$ is a weight simple manifold,
$\hat\nu_i$ is determined by $\nu_i$. Indeed, $T$ acts on $\nu_0$ via the character $\chi\:T\to S^1$
composed with complex multiplication on the fibers. 
The same holds for $\nu_1$, replacing $\chi$ by $\chi^\mun$.
\end{Remark}

Let $(M,M_0,M_1)$ be a \simple\ $T$-manifold with residual moment map
$\bar\Phi\:M\to\bbr$. Let $\ell>0$ defined by $\{\ell\}=\bar\Phi(M_1)$.
Define 
\begin{eqnarray}\label{eq:Vi}
V_0=\bar\Phi^\mun([0,\ell/2])  \ \hbox{ and } \ 
V_1=\bar\Phi^\mun([\ell/2,\ell]) \, .
\end{eqnarray}

\begin{Lemma}\label{voistub}
For $i=0$ and $1$, the subspace $V_i$ of \eqref{eq:Vi} is a $T$-invariant (closed) tubular neighborhood of $M_i$ in $M$.
\end{Lemma}

\preu
We prove this for the case $i=0$, and mention the necessary adaptations
to complete the case $i=1$. The proof introduces
techniques which are useful in subsequent sections (see \remref{R.voistub}
for the idea of a more direct argument). 
Passing to the residual action, we suppose that $T=S^1$.

Choose an $S^1$-invariant almost complex structure $J$ on $M$. This makes $\nu_0$ 
an $S^1$-equivariant Hermitian bundle with structure group $U(r_0)$. 
We denote by $E(\nu_0)$ its total space and by $p\:E(\nu_0)\to M_0$ the bundle projection.
Denote by $S(\nu_0)\subset E(\nu_0)$ the 
associated unit sphere bundle.
For $\varepsilon>0$, let $D_\varepsilon(\nu_0)\subset E(\nu_0)$ the disk bundle formed by the 
elements of $E(\nu_0)$ of norm $\leq\varepsilon$. An element of $D_\varepsilon(\nu_0)$
may be written under the form $rz$, with $z\in S(\nu_0)$ and $r\in[0,\varepsilon]$, with the identification $0z=p(z)$. 

As $\nu_0$ is a Hermitian bundle, each fiber of $E(\nu_0)$ carries a symplectic form,
isomorphic to the standard form on $\bbc^{r_0}$ via a trivialization. The orthogonal sum
with the symplectic form on $M_0$ provides a symplectic form $\omega^0$ on $E(\nu_0)$.
The same construction works for the almost complex structure and the Riemannian metric,
so there is a compatible triple $(\omega^0,J^0,\llangle{\,}{}^0)$ over $E(\nu_0)$,
extending the given one over $M_0$.

Let $b\:D_\varepsilon(\nu_0)\to M$ be the $S^1$-equivariant tubular neighborhood embedding given
by the exponential with respect to the Riemannian metric $\llangle{\,}{}$,
for $\varepsilon>0$ small enough. 
The two symplectic forms $\omega^0$ and $b^*\omega$ coincide on $M_0$.
By \cite[Lemma~3.14]{MDS}, there is a tubular neighborhood embedding 
$h\:D_{\varepsilon'}(\nu_0)\to D_\varepsilon(\nu_0)$ such that $h^*b^*\omega=\omega^0$. 
Based on Moser's argument, the construction of $h$ can be made $S^1$-invariant
(see, e.g. \cite[Remark~II.1.13]{Au}). Thus, replacing $b$ with $b\pcirc h$ and $\varepsilon$ with
$\varepsilon'$ if necessary, we may assume that $b^*\omega=\omega^0$. 
Pushing the triple $(\omega^0,J^0,\llangle{\,}{}^0)$ down to $M$ via $b$,  we get a 
compatible triple $(\omega,J^0,\llangle{\,}{}^0)$ near $M_0$.

Choose a smooth function $\delta_0\:[0,\ell]\to [0,1]$ which is equal to $0$ near $0$ 
and so that the support of $(1-\delta_0)\pcirc\bar\Phi$ is contained in the interior of $b(D_\varepsilon(\nu_0))$. 
Recall that the space $\calj(b(D_\varepsilon(\nu_0)),\omega)$ of $S^1$-invariant $\omega$-compatible almost complex structures on $b(D_\varepsilon(\nu_0))$ is contractible. The standard 
proof of this, for example in
\cite[Propositions~4.1 and 2.49]{MDS}, actually provides a path 
$J^s$ ($s\in [0,1]$) from $J^0$ to $J^1=J$. The formula 
$$
J'_x=J^{\delta_0\pcirc\bar\Phi(x)}_x \in {\rm Aut}_\bbr\, T_xM  
$$
makes sense for all $x\in M$ and provides a $\omega$-compatible almost complex structure
on $M$. We say that $J'$ is obtained by \textit{straightening $J$ around $M_0$}, using 
the \textit{straightening function} $\delta_0$. The almost complex structure $J'$
determines a Riemannian metric $\llangle{\,}{}'$ on $M$, and hence we have an
$S^1$-invariant compatible triple $(\omega,J',\llangle{\,}{}')$ on $M$.

Let us consider the gradient vector field $\grad\bar\Phi$ for the metric $\llangle{\,}{}'$.
This vector field depends only on $J'$, since 
${\rm grad\,}\bar\Phi=J'X$, where $X$ is the fundamental vector field of the Hamiltonian 
residual circle action.
A $J'$-gradient line is the closure of a trajectory of $\grad\bar\Phi$. 

Suppose that $M$ is a weight simple manifold. We claim that for each vector $z\in S(\nu_0)$,
there is a unique $J'$-gradient line $\Gamma_z$ that is tangent to $z$ and 
that hits $M_0$ at a point $p(z)$.
This process parametrizes the gradient lines by $S(\nu_0)$. To see this,
we transport ourselves into $D_\varepsilon(\nu_0)$ via $b$. If $M$ is a weight manifold,
the restriction of the moment map $\bar\Phi\pcirc b$ on each fiber is just the norm square,
whose level surfaces of $\bar\Phi\pcirc b$ are round spheres
and the $J'$-gradient lines are the radial lines to the zero sections. Checking this
also makes it clear that the equation
\begin{equation}\label{E-grad}
\beta_0(rz)=\Gamma_z\cap\bar\Phi^\mun\left(\frac{\ell \,r^2}{2}\right)
\end{equation}
defines a map $\beta_0\:D_1(\nu_0)\to M$ which is an $S^1$-equivariant smooth embedding
with image $V_0$. This completes the proof of \lemref{voistub} for $i=0$ when $M$ is a weight
simple manifold. The case $i=1$ is analogous.
We reverse the orientation of the gradient lines, and for 
$rz\in [0,\sqrt{\ell}]\times\bbd_1$, we define $\beta_1(rz)$ to be
the point $y\in\Gamma_z$ such that $\bar\Phi(y)=\frac{\ell-r^2}{2}$.

Finally, when $M$ is not a weight manifold, the level surfaces of $\bar\Phi\pcirc b$ are ellipsoids 
and the above process does not work: it requires that the Hessian of $\bar\Phi\pcirc b$
be proportional to the metric $\llangle{\,}{}^0$. To get around this difficulty,
we precompose $b$ with an automorphism of $\nu_0$ which transforms the ellipsoids into
round spheres. We use this new tubular neighborhood $b''\:D_{\varepsilon''}\to M$
to transport the metric $\llangle{\,}{}^0$ on a neighborhood of $M_0$ in $M$,
providing a Riemannian metric $\llangle{\,}{}''$ on this neighborhood. This metric 
may be mixed with $\llangle{\,}{}$ using a function like $\delta$ 
to obtain an $S^1$-invariant Riemannian metric $\llangle{\,}{}^-$ on $M$.
Then Equation~\eqref{E-grad} together with the metric $\llangle{\,}{}^-$ 
provides an $S^1$-invariant smooth 
tubular neighborhood embedding with image $V_0$. Note that the metric
$\llangle{\,}{}^-$ is no longer compatible with the symplectic form, but this is not necessary
for the proof of  \lemref{voistub}. 
\cqfd

\begin{rem}\label{R.voistub}
The above proof of \lemref{voistub} was designated to introduce techniques 
useful in subsequent sections.
For a more direct proof, recall that the Morse Lemma provides an embedding $\psi\:D_1(\nu_0)\to M$
with image a tubular neighborhood $\cald$ of $M_0$, such that each gradient line of $\bar\Phi$  
intersects the boundary of $\cald$ transversally in one point. This enables us to construct a
diffeomorphism $\beta_0\:D_1(\nu_0)\to V_0$ as in \eqref{E-grad}.
Thus, $V_0$ is a tubular neighborhood of $M_0$ (note that $V_0$ is $T^1$-invariant by definition). 
\end{rem}

\section{Cohomology constraints}\label{se:coh}

In this paper, $H^*(\cdot)$ denotes the cohomology ring of a space with rational coefficients.
Recall that, for $(M,M_0,M_1)$ a simple Hamiltonian manifold,
$2m_i$ and $2m$ are the dimensions of $M_i$ and $M$ respectively, 
and that $2r_i={\rm codim\,}M_i$.

\begin{Proposition}\label{HKexseq}
Let $(M,M_0,M_1)$ be a \simple\ manifold. 
Then
for $i,j\in\{ 0,1\}$ and $i\neq j$, there are short exact sequences
\begin{equation}\label{HKexseq1}
0 \to H^{*-2r_j}(M_j)\to H^*(M)\to H^*(M_i)\to 0
\end{equation}
and
\begin{equation}\label{HKexseq2}
0 \to H_{2m-*}(M_j)\to H^*(M)\to H^*(M_i)\to 0 ,
\end{equation}
where the right hand homomorphisms are induced by inclusion.
\end{Proposition}

\begin{rem}
This is related to the results in  \cite[\S\,3]{HK}.  Here we do not need to assume that the 
cohomology of $M_0$ and $M_1$ is concentrated in even degrees because the moment
map provides a perfect Morse-Bott function that allows us to deduce the result.
\end{rem}

\preu 
Let $V_j$ be the tubular neighborhood near $M_j$ given by \lemref{voistub} that
satisfies $V_i=M-{\rm int\,}V_j$. 
We first note that the cohomology exact sequence of the pair $(M,V_i)$ 
splits into short exact sequences
\begin{equation}\label{HKexseq2-eq10}
0 \to H^{*}(M,V_i)\to H^*(M)\to H^*(V_i)\to 0 \, .
\end{equation}
This is related to the fact that the residual moment map 
is a perfect Morse-Bott function. A proof of~\eqref{HKexseq2-eq10}
for the $T$-equivariant cohomology is given in \cite[Proposition~2.1]{TW}.
Exactness of \eqref{HKexseq2-eq10} then follows because  $M_i$ is $T$-fixed,
so the map 
$H^*_T(M_i)\to H^*(M_i)$ is onto.
By excision of ${\rm int\,}V_i$ and the Thom isomorphism,
\begin{equation}\label{HKexseq2-eq20}
H^*(M,V_i)\approx H^*(V_j,\partial V_j)\approx H^{*-2r_j}(M_j) \, .
\end{equation} 
Then \eqref{HKexseq2-eq10} and \eqref{HKexseq2-eq20}
give exactness of Sequence~\eqref{HKexseq1}.

Next,
Poincar\'e duality for $V_j$ implies that 
\begin{equation}\label{HKexseq2-eq30}
H^*(M,V_i)\approx H^*(V_j,\partial V_j)\approx H_{2m-*}(V_j)\approx H_{2m-*}(M_j) \, .
\end{equation} 
Thus \eqref{HKexseq2-eq10} and \eqref{HKexseq2-eq30}
imply exactness of Sequence~\eqref{HKexseq2}.
\cqfd

Let $P, P_{i}\in\bbz[t]$ be the Poincar\'e polynomials of $M$ and $M_i$.

\begin{Corollary}\label{HKexseq-cor}
The Poincar\'e polynomial $P_0$ together with $r_0$ and $r_1$ determine
both $P_1$ and $P$ by the following equations:
\begin{equation}\label{HKeq2}
\left\{\begin{array}{rcl}
(1-t^{2r_1})P_1 &=&  (1-t^{2r_0})P_0        \\[2mm]
(1-t^{2r_1})P &=& (1-t^{2(r_0+r_1)})P_0
\end{array}\right. .
\end{equation}
\end{Corollary}

\preu
Sequences~\eqref{HKexseq1} for $i=0$ and $i=1$ immediately give the following
equations
\begin{equation}\label{HKeq2-preu}
\left\{\begin{array}{rcrrrr}
P &=& P_0 &+& t^{2r_1}P_1  \\[2mm]
P &=& t^{2r_0}P_0 &+& P_1
\end{array}\right. ,
\end{equation}
from which we may deduce the equations of \corref{HKexseq-cor}. Note that Equations~\eqref{HKeq2-preu}
are just the Morse-Bott equalities for the residual moment map and its opposite. 
\cqfd

\begin{Corollary}\label{HKexseq-cor05}
Let $(M,M_0,M_1)$ be a simple Hamiltonian manifold
with $r_1=2$. Then there are additive isomorphisms
$$
H^*(M_1)\approx_{{\rm add}}H^*(M_0)\otimes H^*(\bbc P^{r_0-1})
\ \hbox{ and } \
H^*(M)\approx_{{\rm add}}H^*(M_0)\otimes H^*(\bbc P^{r_0}).
$$
\end{Corollary}

\preu
Suppose that $M$ is obtained by a symplectic cut of the trivial
bundle $M_0\times\bbc^{r_0}$. Then $M_1=M_0\times\bbc P^{r_0-1}$
and $M=M_0\times\bbc P^{r_0}$, which proves the lemma in this case.
The general case follows from \corref{HKexseq-cor}.
\cqfd

\begin{Remark} \rm
It is not true that $P_0$ together with $r_1$ determines
the cohomology ring $H^*(M)$. For instance, for the symplectic cut of a weight bundle
$\nu_0$ over $M_0$ given in Example~\ref{eg:sympcut}, the ring structure on $H^*(M)$ depends
on the bundle $\nu_0$. For $M_0=S^2$ and $r_0=1$, $M$ is diffeomorphic to $S^2\times S^2$ if
$c_1(\nu_0)$ is even and to $\bbc P^2\sharp\,\overline{\bbc P}^2$ if $c_1(\nu_0)$ is odd. 
\end{Remark}

The first equation in~\eqref{HKeq2} immediately implies the following corollary.

\begin{Corollary}\label{HKexseq-cor10}
If $r_0=r_1$, the Poincar\'e polynomials of $M_0$ and $M_1$ are identical:  $P_0=P_1$. \cqfd
\end{Corollary}

The following proposition appears as a special case of the first centered equation in \cite{HT}.  In the case of a simple Hamiltonian manifold, their inequality is precisely this one.

\begin{Proposition}\label{m1m2m}
Let $(M,M_0,M_1)$ be a simple Hamiltonian manifold. Then
$$
m\leq m_0+m_1+1
$$
\end{Proposition}

\preu Suppose that $2r_1> 2m_0+2$. The first equation of
\eqref{HKeq2-preu} then implies that $H^{2m_0+2}(M)=0$, which is impossible
as $M$ is a compact symplectic manifold of dimension $\geq 2m_0+2$.
Hence, $2r_1\leq 2m_0+2$,
which implies that $2m\leq 2m_0+ 2m_1+2$.
\cqfd

\begin{Lemma}\label{L.H1}
Let $(M,M_0,M_1)$ be a simple Hamiltonian manifold. Then
$$
H^1(M_0)\approx H^1(M) \approx H^1(M_1),
$$
these isomorphisms being induced by the inclusions $M_i\subset M$.
\end{Lemma}

\preu
As $m_i\geq 1$,
the abstract isomorphisms come from Equations~\eqref{HKeq2-preu}.
By \proref{HKexseq}, inclusions $M_i\subset M$ induce
surjective homomorphisms, which are then isomorphisms.
\cqfd

\begin{Proposition}\label{P.oddcoh}
For a simple Hamiltonian manifold $(M,M_0,M_1)$, the following conditions
are equivalent.
\renewcommand{\labelenumi}{(\alph{enumi})}
\begin{enumerate}
\item $H^{odd}(M_0)=0$.
\item $H^{odd}(M_1)=0$.
\item $H^{odd}(M)=0$.
\end{enumerate}
\end{Proposition}

\preu
By the first equation of  \eqref{HKeq2},
Conditions (a) and (b) are equivalent.
By Equation~\eqref{HKeq2-preu}, (c) is equivalent to
(a) and (b) together.
\cqfd

\begin{Example}\rm
Suppose that $M$ has the cohomology of $\bbc P^n$. Then $M_0$ and
$M_1$ have the cohomology ring of a complex projective space.
Indeed, their cohomology groups vanish in odd degree by
\proref{P.oddcoh}. Also, their Betti numbers are $\leq 1$ by
\proref{HKexseq} and they are symplectic manifolds.
The first equation of \eqref{HKeq2-preu}.
implies that $m_0 + m_1+1=m$, as in Example~\ref{edgeexfirst}
(For $M_1=pt$, this is a result of \cite[Theorem~1]{HR}). 
\end{Example}

\begin{Remark}\label{R.HUI-Tolman}\rm
The extreme case in \proref{m1m2m}, i.e. $m= m_0+m_1+1$, is studied in
\cite{T}, \cite{HT} and \cite{LOS}. Much stronger restrictions than what we prove in this section hold
in that special case.  In that context, the ring $H^{*}(M;\bbz)$ must be isomorphic
either to  $H^*(\bbc P ^m)$ or  to $H^*(\tilde G_2(\bbr^{m+2}))$,
and $M$ is necessarily simply connected. 
Moreover, $M_1$ and $M_2$ each have the homotopy type of a complex projective space.
\end{Remark}

\section{Diffeomorphism invariants}\label{S.Diffeo}

Let $M_0^a$ and $M_1^a$ be fixed compact smooth manifolds (the exponent $a$ stands for 
\textbf{abstract}). We also fix two Hermitian vector bundles $\nu_i^a\: \bbe_i\to M_i^a$ of complex rank $r_i$. The isomorphism class $[\nu_i^a]$ of the \textbf{abstract normal bundle} may be considered as an element of $[M_i^a,BU(r_i)]$; we write 
$$[\nu^a]=([\nu_0^a],[\nu_1^a])\in [M_0^a,BU(r_0)]\times[M_1^a,BU(r_1)].$$

\begin{definition}\label{def:normal-simple}
A \textbf{$[\nu^a]$-simple Hamiltonian} $T$-manifold consists of a weight simple Hamiltonian $T$-manifold 
$(M,M_0,M_1)$ together with diffeomorphisms $\alpha_i\:M_i^a\fl{\approx} M_i$ for $i=0,1$,
such that $\alpha_i^*[\nu_i]=[\nu_i^a]$. Here, $\nu_i=TM|_{M_i}/TM_i$ is called the \textbf{concrete} normal bundle 
to $M_i$ in $M$.  It can be endowed with a $U(r_i)$-structure group via the choice of 
an almost complex structure $J\in \calj(M,\omega)$. 
\end{definition}

The isomorphism class $[\nu_i]$ is well-defined (see the discussion before \remref{R.Char}).
Two such objects $((M,M_0,M_1),\alpha_i)$ and $((M',M_0',M_1'),\alpha_i')$ are 
considered equivalent if there is a $T$-equivariant symplectomorphism 
$h\:M\to M'$ such that $h\pcirc\alpha_i=\alpha_i'$. The set of equivalence classes
of $[\nu^a]$-simple Hamiltonian $T$-manifolds
is denoted $\calh([\nu^a])$.

The first invariant associated to a class $\calm\in\calh([\nu^a])$ is the 
character $\chi(\calm)\in \hat T$
defined in \lemref{character}. Note that, since $M\neq M^T$, the map $\chi:T\to S^1$ is surjective. 
As we are dealing with weight manifolds, the residual action is semi-free, with
residual moment map: $\bar\Phi\:M\to [0,\ell]$,
that sends $M_0$ to $0$. The number $\ell=\ell(\calm)>0$ is another invariant of the class  
$\calm\in\calh([\nu^a])$, called the \textbf{$T$-size} of $\calm$.

Note that $\nu_i^a$ and the character $\chi$ determine unique $T$-equivariant weight bundles,
as discussed in \remref{R.Char}. Thus, $\nu_i^a$ is $T$-equivariantly isomorphic to the concrete normal bundle
$\nu_i$ of a representative of $\calh([\nu^a])$.
Associated to the abstract normal Hermitian bundle $\nu_i^a$, we have the following.

\begin{definition}\label{def:laundry-list}
For the bundle $\nu_i^a$, denote the total space $\bbe_i$ with its bundle projection $p_i\:\bbe_i\to M_i^a$.
This has associated bundles and structure groups:
\begin{enumerate}
\item[\ref{def:laundry-list}.1] the \textbf{abstract sphere bundle} $\bbs_i\to M_i^a$ (fiber $S^{2r_i-1}$),
where $$\bbs_i=\{z\in\bbe_i\mid \, |z|=1\}.$$

\item[\ref{def:laundry-list}.2] the \textbf{abstract disk bundle} $\bbd_i\to M_i^a$ (fiber the unit disk in $\bbc^{r_i}$),
where $\bbd_i=\{z\in\bbe_i\mid \, |z|\leq1\}$.
We also consider the disk bundle $\bbd_{i,\varepsilon}=\{z\in\bbe_i\mid \, |z|\leq\varepsilon\}$.

\item[\ref{def:laundry-list}.3] the \textbf{abstract projective bundle}  $\bbp_i\to M_i^a$ (fiber $\bbc P^{2r_i-1}$),
where $\bbp_i=\bbs_i/S^1$. The projection
$\eta_i\:\bbs_i\to\bbp_i$ is a principal $S^1$-bundle with Euler class $e(\eta_i)\in H^2(\bbp_i;\bbz)$.

\item[\ref{def:laundry-list}.4] the \textbf{extended gauge group} $\hat\calg(\nu_i^a)$, defined by pairs of isomorphisms
that fit into commutative diagrams
$$
\begin{array}{c}{\xymatrix@C-3pt@M+2pt@R-4pt{%
\bbe_i \ar[d]
\ar[r]^(0.50){g}  &
\bbe_i \ar[d]  \\
%%%%%% ROW 2
M_i^a \ar[r]^(0.50){\bar g}  &
M_i^a
}}\end{array},
$$
where $g$ is smooth and its restriction to each fiber is an isometry. 
Those isomorphisms with $\bar g= {\rm id}$
form the usual \textbf{gauge group} $\calg(\nu_i^a)$.
There is thus an exact sequence
\begin{equation}\label{seqgbarg}
1 \to \calg(\nu_i^a) \to \hat\calg(\nu_i^a) \to {\rm Diff\,}(M_i^a,[\nu_i^a]) \to 1,
\end{equation}
where ${\rm Diff\,}(M_i^a,[\nu_i]^a)$ denotes the group of diffeomorphisms $h\:M_i^a\to M_i^a$
that satisfy $h^*[\nu_i^a]=[\nu_i^a]$.
The group $\hat\calg(\nu_i^a)$ acts naturally on each of the above associated bundles.

\item[\ref{def:laundry-list}.5] the \textbf{extended gauge group} $\hat\calg(\eta_i)$, defined by pairs of isomorphisms 
that fit into commutative diagrams
$$
\begin{array}{c}{\xymatrix@C-3pt@M+2pt@R-4pt{%
\bbs_i \ar[d]
\ar[r]^(0.50){g}  &
\bbs_i \ar[d]  \\
%%%%%% ROW 2
\bbp_i \ar[r]^(0.50){\bar g}  &
\bbp_i
}}\end{array}
$$
such that $g$ is smooth and $S^1$-equivariant. 
Those isomorphisms with $\bar g= {\rm id}$
form the usual \textbf{gauge group} $\calg(\eta_i)$. 
\end{enumerate}
\end{definition}

The $T$-action on $\nu_i^a$ induces a $T$-action on all the abstract sphere and disk bundles 
which commutes with the actions of the extended gauge groups.

Let $((M,M_0,M_1),\alpha_i)$ represent an element of $\calh([\nu^a])$. 
Choose a compatible almost complex structure $J$ on $M$, and consider its associated 
Riemannian metric. 
As discussed above, this endows the concrete normal bundle 
$\nu_i=TM|_{M_i}/TM_i$ with a $T$-invariant
Hermitian structure, making it isometric to the orthogonal complement of $TM_i$ in $TM$.
Choose Hermitian vector bundle isomorphisms $\gamma_i\:\bbe_i\to E(\nu_i)$ 
covering $\alpha_i$. 
These induce isomorphisms on the associated bundles:
$\gamma_i\:\bbs_1\to S(\nu_i)$ and so forth.
We also get a tubular neighborhood embedding $b\:\bbd_{i,\varepsilon}\to M$
of $M_i$ in $M$. We now proceed as in the proof of \lemref{voistub}. 
We may use the embedding $b$ to
straighten the Riemannian metric around $M_i$, using straightening functions 
$\delta_i\:[0,\ell]\to [0,1]$. 
The gradient lines for the moment map $\Phi$ and the straightened metric 
provide a $T$-equivariant smooth embedding
$\beta_0\:\bbd_0\to M$ by 
$$
\beta_0(rz)=\Gamma_{\gamma_0(z)}\cap \bar\Phi^\mun\left(\frac{\ell \,r^2}{2}\right)
$$ 
where $\Gamma_{\gamma_0(z)}$ is the unique gradient line starting from $p_0(z)$ in the direction
of $\gamma_0(z)$. 
The $T$-equivariant embedding $\beta_1\:\bbd_1\to M$ is defined symmetrically.
The image of $\beta_0$ and $\beta_1$ are the $T$-invariant tubular neighborhoods
$V_0=\phi^\mun([0,\ell/2])$ and $V_1=\phi^\mun([\ell/2,\ell])$. 

The map 
\begin{equation}\label{difpsi}
\psi=\beta_0^\mun\pcirc\beta_1\:\bbs_1\to\bbs_0
\end{equation}
is a diffeomorphism which anti-commutes with the
$S^1$-action.
Let $\cale(\nu^a)$ be the space of such diffeomorphisms $\psi\:\bbs_1\to\bbs_0$. 
Observe that $\psi$ descends to a diffeomorphism $\bar\psi\:\bbp_1\to\bbp_0$.
By pre-composition, the extended gauge group $\hat\calg(\nu_1^a)$ acts on the right on $\cale(\nu^a)$
and, by post-composition, $\hat\calg(\nu_0^a)$ acts on the left on $\cale(\nu^a)$.
These two actions commute and descend to the isotopy classes,
giving actions of $\pi_0(\hat\calg(\nu_1^a))$ and $\pi_0(\hat\calg(\nu_0^a))$
on $\pi_0(\cale(\nu^a))$. We can restrict these actions to the usual gauge groups.
Define the set $\cale([\nu^a])$ by 
\begin{equation}\label{dcoset-1}
\cale([\nu^a])=\pi_0(\calg(\nu_0^a)) \big\backslash \pi_0(\cale(\nu^a)) \big/ 
\pi_0(\calg(\nu_1^a)) 
\end{equation}
The notation $\cale([\nu^a])$ makes sense because the 
above double coset
\textit{depends only on} $[\nu^a]$. More precisely,
let $\nu'=(\nu_0',\nu_1')$ and  $\nu''=(\nu_0'',\nu_1'')$ be two representatives
of $[\nu^a]$. Choosing principal bundle isomorphisms $\kappa_i\:E(\nu_i')\to E(\nu_i'')$
produces a bijection $\kappa$ between the double quotient~\eqref{dcoset-1}
for $\nu'$ and $\nu''$. Since we have divided out by the action of the gauge groups, 
the bijection $\kappa$ does not depend on the choice
of the $\kappa_i$'s.

\begin{Lemma}\label{L.psiwdef}
The above construction provides a well-defined map
$$\Psi\:\calh([\nu^a]) \to \cale([\nu^a]) \, .$$
\end{Lemma}

\preu 
Let $((M,M_0,M_1),\alpha_i)$ 
represent a class $\calm\in\calh([\nu^a])$. 
The definition of the diffeomorphism $\psi$ of~\eqref{difpsi}
involves three choices:
\begin{enumerate}
\renewcommand{\labelenumi}{(\alph{enumi})}
\item the compatible almost complex structure $J$ on $M$;
\item the $U(r_i)$-isomorphism $\gamma_i:\bbe_i\to E(\nu_i)$; and
\item the straightening functions $\delta_i\:[0,\ell]\to [0,1]$.
\end{enumerate}
Once the choices (a) and (b) have been made, the straightening functions $\delta_0$ and $\delta_1$
belong to convex spaces, so their choice does not change $\psi$ in $\pi_0(\cale(\nu))$.
If we choose instead $\tilde\gamma_i\:\bbe_i\to E(\nu_i)$ for (b), then 
$\tilde\gamma_i=g_i\pcirc\gamma_i$ with $g_i\in\calg(\nu_i)$. Hence,
$\tilde\psi= g_0\pcirc\psi \pcirc g_1$, proving that $\tilde\psi$ and $\psi$ represent the
same class in $\cale([\nu])$. Finally, the choice of (a) does not change the class
since compatible almost complex structures on $M$ form a contractible space.

Now, let $((\bar M,\bar M_0,\bar M_1),\bar \alpha_i)$ 
be another representative of $\calm$.
Let $h\:M\to \bar M$ be a $T$-equivariant symplectomorphism realizing the
equivalence. Choose a compatible almost complex 
structure $\bar J$ on $\bar M$. Then $J=Th^\mun\pcirc\bar J\pcirc Th$ is a compatible almost complex 
structure on $M$, which may be used, together with  the above Hermitian bundle isomorphisms
$\gamma_i$ to get a representative $\psi$ of $\Psi(\calm)$. 
The construction is transported via $h$ to $\bar M$, using $\bar J$, setting
$\bar\gamma_i=Th\pcirc\gamma_i$, 
and using the same straightening function. We thus get
embeddings $\bar\beta_i=h\pcirc\beta_i\:\bbd_i\to M'$ which can be used to define
$\bar\psi\:\bbs_1\to\bbs_0$, which then satisfies
$$
\bar\psi = \bar\beta_0^\mun\pcirc\bar\beta_1 =  \beta_0^\mun\pcirc h^\mun\pcirc h \pcirc\beta_1 = \psi \,.
\cqfd
$$

Let us consider the following quotients of the set $\cale([\nu^a])$:
\begin{equation}\label{dcoset-2}
\cale^{1}([\nu^a])=\pi_0(\calg(\nu_0^a)) \big\backslash \pi_0(\cale(\nu^a)) \big/ 
\pi_0(\hat\calg(\nu_1^a)) 
\end{equation}
and
\begin{equation}\label{dcoset-3}
\cale^{01}([\nu^a])=\pi_0(\hat\calg(\nu_0^a)) \big\backslash \pi_0(\cale(\nu^a)) \big/ 
\pi_0(\hat\calg(\nu_1^a)) 
\end{equation}
The compositions of the map $\Psi\:\calh(\nu^a])\to\cale([\nu^a])$ with the projections
onto $\cale^{1}([\nu^a])$ and $\cale^{01}([\nu^a])$ are denoted by $\Psi^{1}$ and $\Psi^{01}$.

\begin{Theorem}\label{TDIFF}
Let $((M,M_0,M_1),\alpha_i)$ and $((M',M_0',M_1'),\alpha_i')$ represent classes 
$\calm$ and $\calm'$ in $\calh([\nu^a])$, with $\chi(\calm)=\chi(\calm')$.
Denote the reduced moment maps by $\bar\Phi\:M\to [0,\ell]$ and $\bar\Phi'\:M\to [0,\ell']$,
where $\ell$ and $\ell'$ are the $T$-sizes of $\calm$ and $\calm'$.
\begin{enumerate}
\renewcommand{\labelenumi}{(\alph{enumi})}
\item If we have $\Psi(\calm)=\Psi(\calm')$, then there is a $T$-equivariant diffeomorphism
$h\:M\to M'$ satisfying 
\begin{equation}\label{TDIFF-eq}
\bar\Phi'\pcirc h = \frac{\ell'}{\ell\hskip 1.2mm} \,\bar\Phi 
\end{equation}
and such that $h\pcirc\alpha_i=\alpha_i'$ for $i=0,1$. 
\item  If we have  $\Psi^{1}(\calm)=\Psi^{1}(\calm')$, then there is a $T$-equivariant diffeomorphism
$h\:M\to M'$ satisfying~\eqref{TDIFF-eq}
and such that $h\pcirc\alpha_0=\alpha_0'$. 

\item If $\Psi^{01}(\calm)=\Psi^{01}(\calm')$, then there is a $T$-equivariant diffeomorphism
$h\:M\to M'$ satisfying~\eqref{TDIFF-eq}. 
\end{enumerate}
\end{Theorem}

Equation~\eqref{TDIFF-eq} means that $\Phi'\pcirc h=\sigma\pcirc \Phi$ where $\sigma$ is 
an affine isomorphism of $\algt^*$ of ratio $\ell'/\ell$.

\sk{2}\preu
For Part (a), choose $(J,\gamma_i)$ and $(J',\gamma_i')$ as above,
getting $T$-equivariant embeddings $\beta_i$ and $\beta_i'$ and $\psi,\psi'\in\cale(\nu)$.
The condition $\Psi(\calm)=\Psi(\calm')$ implies
an equation in $\pi_0(\cale(\nu))$ of the form $[\psi']=g_1[\psi]g_0$ with $g_i\in\calg_i$.
Changing $\gamma_i$ into $\gamma_i\pcirc g_i^\mun$, we get that $[\psi]=[\psi']$
in $\pi_0(\cale(\nu))$. Now, the embeddings $\beta_i$ produce a
$T$-equivariant diffeomorphism $N_\psi=\bbd_0\cup_\psi\bbd_1 \fl{q} M$ 
extending $\alpha_0$ and $\alpha_1$. 
In the same way, the embeddings $\beta_i'$ produce a
$T$-equivariant diffeomorphism 
$$N_{\psi'}=\bbd_0\cup_{\psi'}\bbd_1\fl{q'} M'$$
extending $\alpha_0'$ and $\alpha_1'$.
As $[\psi]=[\psi']$,
there is a smooth $T$-equivariant isotopy $$b\:\bbs_0\times [1/2,1]\to \bbs_0\times [1/2,1],$$ preserving
the projection onto $[1/2,1]$,  such that $b(z,t)=(z,t)$ for $t$ near $1/2$,
and $b(z,t)=(\psi'\pcirc\psi^\mun(z),t)$ for $t$ near $1$. This isotopy extends, by the identity 
near the null-section, to a $T$-equivariant diffeomorphism $b\:\bbd_0\to\bbd_0$. Now, $b$ together with the identity on $\bbd_1$ gives a $T$-equivariant diffeomorphism $B\:N_\phi\fl{\approx} N_{\phi'}$.
Finally, observe that the level sets of the maps $\bar\Psi\pcirc q$ and $\bar\Psi\pcirc q'$ are
the manifolds $|z|=constant$ in $\bbd_i$. These level sets are preserved by the diffeomorphism 
$B$. 
By the definition of the embeddings $\beta_i$ and $\beta_i'$, this proves Equation~\eqref{TDIFF-eq}
and completes the proof of (a).
Parts (b) and (c) are proven in the same way, but the elements $g_i$ that occur in the
above argument are now in $\hat\calg_i$ instead of $\calg_i$.
\cqfd

In order to get applications of \thref{TDIFF}, we
now provide a different description of $\cale([\nu])$ and its quotients.
Choose an element $h\in\cale([\nu^a])$, if $\cale([\nu^a])$ is non-empty. 
Then any $\tilde h\in\cale([\nu^a])$
is of the form $\tilde h = h \pcirc (h^\mun\pcirc\tilde h)$ and
$h^\mun\pcirc\tilde h\in\hat\calg(\eta_1)$. Hence, the map
$g\mapsto h\pcirc g$ provides a bijection from $\hat\calg(\eta_1)$ onto  $\cale([\nu])$.
Now, there is an injection $\hat\calg(\nu_0)\hookrightarrow \cale([\nu])$ given by
$\gamma\mapsto \gamma\pcirc h$. Composed with the above bijection
$\hat\calg(\eta_1)\fl{\approx}\cale([\nu^a])$ gives an injective
homomorphism $$\hat\calg(\nu_0^a)\to\hat\calg(\eta_1)$$ defined by
$\gamma\mapsto h\pcirc \gamma\pcirc h^\mun$.
We have proven the following proposition.

\begin{Proposition}\label{compucale}
If $\cale([\nu]^a)$ is not empty, the choice of $h\in\cale(\nu^a)$ provides
bijections
$$\cale([\nu^a]) \fl{\approx}
\pi_0(\calg(\nu_0^a))\big\backslash\pi_0(\hat\calg(\eta_1))\big/\pi_0(\calg(\nu_1^a)) \, ,$$ 
$$\cale^{1}([\nu]^a) \fl{\approx}
\pi_0(\calg(\nu_0^a))\big\backslash\pi_0(\hat\calg(\eta_1))\big/\pi_0(\hat\calg(\nu_1^a)) 
$$
and 
$$\cale^{01}([\nu]^a) \fl{\approx}
\pi_0(\hat\calg(\nu_0^a))\big\backslash\pi_0(\hat\calg(\eta_1))\big/\pi_0(\hat\calg(\nu_1^a)) \, ,
$$
where the inclusion $\hat\calg(\nu_0^a)\hookrightarrow\hat\calg(\eta_1^a)$ is given 
by $\gamma\mapsto h^\mun\pcirc\gamma\pcirc h$. 
\cqfd
\end{Proposition}

\section{The case $r_1=1$}\label{se:case}

The results of this section follow from the following proposition.

\begin{Proposition}\label{P.r1=1cale}
Let $M_i^a$ be compact smooth manifolds  for $i=0,1$. Let 
$$[\nu^a]=([\nu_0^a],[\nu_1^a])\in [M_0^a,BU(r_0)]\times[M_1^a,BU(r_1)].$$
Suppose that $r_1=1$. Then $\cale^1([\nu^a])$ is either empty or contains a single element.
\end{Proposition}

\preu
If $\cale^1([\nu^a])$ is not empty, then it is, by \proref{compucale}, in bijection
with 
$$\pi_0(\calg(\nu_0^a))\big\backslash\pi_0(\hat\calg(\eta_1))\big/\pi_0(\hat\calg(\nu_1^a)).$$
As $r_1=1$, $\nu_1^a$ is isomorphic to the complex line bundle associated to $\eta_1$.
Hence, $\hat\calg(\nu_1^a)=\hat\calg(\eta_1)$ which implies that $\cale^1([\nu])$ consists
of a single element.
\cqfd

We now provide a criterion to determine, in \proref{P.r1=1cale}, whether 
$\cale^1([\nu^a])$ is non-empty. Let $\eta_i:\bbs_i\to\bbp_i$ be the $S^1$-principal 
bundle associated to $\nu_i^a$. Let $\bbl_i\to\bbp_i$ be the Hermitian line bundle associated to $\eta_i$. Let $\bbl_i^-\to\bbp_i$ be the conjugate line bundle, and denote its isomorphism class
by $[\eta_i^-]$.

\begin{Proposition}\label{P.r1=1cale-crit}
Let $M_i^a$ and $[\nu^a]$ as in \proref{P.r1=1cale}. The set $\cale^1([\nu^a])$ is non-empty
if and only if there exists a diffeomorphism $\kappa\:M_1^a\to\bbp_0$ such that
$\kappa^*[\eta_0^-]=[\nu_1^a]$.
\end{Proposition}

\preu The diffeomorphism $\kappa$ would be covered by a diffeomorphism $\tilde\kappa\:\bbs_1\to\bbs_0$
which anti-commutes with the $S^1$-action. Such a $\tilde\kappa$ defines a class
in $\cale^1([\nu^a])$.

Conversely, a class in $\cale^1([\nu^a])$ is represented by 
a diffeomorphism $h\:\bbs_1\to\bbs_0$ which anti-commutes with the $S^1$-action.
This descends to $\bar h\:\bbp_1\to\bbp_0$ satisfying $\bar h^*[\eta_0^-]=[\eta_1]$. 
As $r_1=1$, there is a bundle isomorphism between $\bbe_1$ and $\bbl(\eta_1)$ (over the identity
of $M_1^a$). Hence, $\kappa=\bar h$ is the desired diffeomorphism.
\cqfd

We now describe in details a basic example.

\begin{Example}\label{E.scut}\rm
Let $N$ be a compact symplectic manifold.
Let $\xi\:E\to N$ be a Hermitian vector bundle of complex rank $r$.
Each fibre of $\xi$ is equipped with a symplectic form coming from the standard symplectic form on $\bbc^{r}$
via a trivialization. Then the symplectic form on $N$ as well as those 
on the fibres of $\xi$ are the restriction of a 
unique symplectic form $\omega$ on $E$. 
The action of $S^1$
by complex multiplication is Hamiltonian, with moment map
$\bar\Phi(z)=\frac{1}{2}||z||^2$. Any $\ell>0$ is a regular value,
so we may take the symplectic cut $\hat P_\ell(\xi)$ of $E$ at $\ell$.
We thus get a simple $S^1$-Hamiltonian manifold
$(\hat P_\ell(\xi),N, P_\ell(\xi))$, where
$P_\ell(\xi)$ is the symplectic reduction of $E$ at $\ell$.
Using a non-trivial character $\chi\:T\to S^1$, we get thus get a weight simple $T$-Hamiltonian manifold
with residual moment map $\bar\phi$. We denote this simple Hamiltonian manifold by 
$\calc_\chi(N,\xi,\ell)$.

Let us define abstract manifolds and normal bundles for 
$\calc_\chi(N,\xi,\ell)$. We can take $M_0^a=N$, $\alpha_0={\rm id}$ and
$\nu_0^a=\xi$. Then there is 
a canonical diffeomorphism $\alpha_1\:\bbp_0\fl{\approx } P_\ell(\xi)$ obtained
by following the real vector lines in $\bbe_0=E$. Hence, together with $\alpha_0$ and $\alpha_1$,
$\calc_\chi(N,\xi,\ell)$ is a $(N,\bbp_0)$-simple
Hamiltonian $T$-manifold. As seen in \proref{P.r1=1cale-crit}, $[\nu_1]=[\eta_0^-]$.

The $T$-embedding $\beta_0$ is induced by the embedding 
$\tilde\beta_0\:\bbd_0\to\bbe_0$ defined by $\tilde\beta_0(rz)=r\sqrt{\ell}\,z$.
Using the identification $\bbs_1=\bbs_0^-$, the elements of $\bbd_1$ may be written under the
form $r\, z$ with $r\in [0,1]$ and $z\in\bbs_0$, with the identification 
$0\,z=0\,z'=p(z)=p(z')$ when the projection of $z$ and $z'$ onto $\bbp_0$ coincide.
The $T$-embedding $\beta_1$ is then induced by the $T$-map 
$\tilde\beta_1\:\bbd_0^-\to\bbe_0$ defined by 
$\tilde\beta_1(rz)=r(\sqrt{\ell}-\sqrt{2\ell})\,\bar z$.
Hence, 
$$
\Psi(\calp)=[{\rm id}]
$$
(the identity from $\bbs_0$ to $\bbs_0^-=\bbs_1$ anti-commuting with the $S^1$-multiplication,
as expected).
\end{Example}

\begin{Theorem}\label{T.r1=1}
Let $(M,M_0,M_1)$ be a simple Hamiltonian $T$-manifold with $T$-size $\ell$,
and associated character $\chi$.
Suppose that $r_1=1$.
Then there exits a $T$-equivariant diffeomorphism 
$$
F\:\calc_\chi(M_0,\nu_0,\ell)\fl{\approx}M
$$
commuting with the residual moment maps and such that $F|_{M_0}= {\rm id}$.
\end{Theorem}

\preu
As $r_1=1$, we know $M$ is a weight simple Hamiltonian manifold.
Define $M_0^a=M_0$ and set $\alpha_0={\rm id}$. Fix an almost complex structure on $M$ compatible with
the symplectic form and let $\nu_0^a$ be the orthogonal complement of $TM_0$ in $TM$
for the associated metric. By \proref{P.r1=1cale-crit},
there exists a diffeomorphism $\alpha_1\:\bbp_0\to M_1$ such that
$\alpha_1^*[\nu_1^a]=[\eta_0^-]$. Hence, 
$((M,M_0,M_1),\alpha_i)$ represents a class in $\calh([\nu])$ for 
$[\nu]=([\nu_0^a],[\eta_0^-])$. So does the simple Hamiltonian manifold
$\calc_\chi(M_0,\nu_0,\ell)$ of \exref{E.scut}, with its own $\alpha_i$'s.
By \thref{TDIFF} and \proref{P.r1=1cale}, this completes the proof of  \thref{T.r1=1}.
\cqfd

\thref{T.r1=1} implies that $M_1$ is diffeomorphic to $\bbp_0$. 
If, in addition $r_0=1$, then $\bbp_0$
is diffeomorphic to $M_0$ and we have the following corollary,
also found in \cite[Lemma~3.2]{GH}.

\begin{Corollary}\label{C-cod2-2}
Let $(M,M_0,M_1)$ be a simple Hamiltonian manifold with 
$r_0=r_1=1$. Then $M_1$ is diffeomorphic to $M_0$.\cqfd
\end{Corollary}

\section{Classification up to $T$-equivariant symplectomorphism}\label{se:symplectomorphism}

The philosophy of this section is slightly different from that in 
Section~\ref{S.Diffeo}.
We fix a single compact smooth manifold $M_0^a$ and
a Hermitian vector bundle $\nu_0^a\: \bbe_0\to M_0^a$ of complex rank $r_0$,
whose isomorphism class is denoted by $[\nu_0^a]\in [M_0^a,BU(r_0)]$.
The associated bundles $\bbs_0\to M_0$ and so forth, as well as $\eta_0$, are defined as in Section~\ref{S.Diffeo}.

\begin{definition}\label{def:simple-Hamiltonian}
A \textbf{$[\nu_0^a]$-simple Hamiltonian} $T$-manifold consists of a weight simple Hamiltonian $T$-manifold $(M,M_0,M_1)$ 
together with a diffeomorphism $\alpha_0\:M_0^a\fl{\approx} M_0$
such that $\alpha_0^*[\nu_0]=[\nu_0^a]$.
\end{definition}

Here, $\nu_0$ is the \textbf{concrete} normal bundle to $M_0$ in $M$, represented
by the orthogonal complement of $TM_0$ in $TM$ for the Riemannian metric associated
to a  $T$-invariant almost complex structure on $M$ compatible with the symplectic form.
In particular, $\omega^a=\alpha_0^*\omega_0$
is a symplectic form on $M_0^a$.
Two such objects $((M,M_0,M_1),\alpha_0)$ and $((M',M_0',M_1'),\alpha_0')$ are
considered as equivalent if there is a $T$-equivariant symplectomorphism
$h\:M\to M'$ such that $h\pcirc\alpha_0=\alpha_0'$.
The following are invariants of an equivalence class:
\begin{itemize}
\item the associated character $\chi$ and the residual action, which is semi-free,
since we are in the case of weight simple Hamiltonian manifolds;
\item  the $T$-size $\ell>0$;
\item the symplectic form $\omega_0^a$ on $M_0^a$; and
\item the codimensions $r_0$ and $r_1$.
\end{itemize}
Fixing $[\nu_0^a]$, $\omega_0^a$, $\ell$ and $r_1$, we get a set of equivalence classes denoted by 
$$\cals^0([\nu_0^a],\omega_0^a,r_1,\ell).$$
We are especially interested in the case $r_1=1$.
By \thref{T.r1=1}, elements of $\cals^0([\nu_0],\omega_0^a,1,\ell)$ are in bijection with classes 
of symplectic
forms on $\calc_\chi(M_0^a,\nu_0^a,\ell)$ coinciding with $\omega_0^a$ on $M_0^a$ and for
which the $T$-action is Hamiltonian. Two
such forms $\omega$ and $\omega'$ are equivalent if there is a self-diffeomorphism
$F$ of $\calc_\chi(M_0^a,\nu_0,\ell)$, commuting with the reduced moment maps, such that $F^*\omega=\omega'$ and $F|_{M_0^a}={\rm id}$.

Let $\Omega^{sym}(M_0^a)$ be the space of symplectic forms on $M_0^a$, with the topology induced
by the $C^\infty$-topology in $\Omega^2(M_0^a)$. Define
$$
\cald((M_0^a,\omega_0^a),[\nu_0^a],\ell)=\bigg\{\omega\:[0,\ell]\to \Omega^{sym}(M_0^a)\,{\bigg|}\,
\omega(0)=\omega_0^a \hbox{ and } [\omega(\lambda)]=[\omega_0^a]+\lambda e(\eta_0) \bigg\} \, ,
$$
where the last equation holds in de Rham cohomology $H^2_{dr}(M_0^a)$.

\begin{Theorem}\label{T.symcod1}
Suppose that $r_1=1$. Then there exists a bijection
$$
\Theta : \cals^0([\nu_0^a],\omega_0^a,1,\ell) \fl{\approx}
\pi_0\big( \cald((M_0^a,\omega_0^a),[\nu_0^a],\ell)\big) \, .
$$
\end{Theorem}

\preu
Let $M=\calc_\chi(M_0^a,\nu_0,\ell)$.
As noted above, a class of $a\in\cals^0([\nu_0^a],\omega_0^a,1,\ell)$ is represented
by a symplectic form $\omega$ on $M$. Observe that there is a diffeomorphism from $M/S^1$ to
$[0,\ell]\times M_0^a$. The first component is given by the residual moment map and the
second one is induced by the projection $\bbe_0\to M_0^a$. Each slice $\{\lambda\}\times M_0$ is
then endowed with a symplectic form $\omega(\lambda)$ given by the symplectic reduction of $\bbe_0$
at $\lambda$. This provides a map $\omega\:[0,\ell]\to \Omega^{sym}(M_0^a)$ with
$\omega(0)=\omega_0^a$. The equation  $[\omega(\lambda)]=[\omega_0^a]+\lambda e(\eta_0)$
holds in $H^2(M_0^a)$ by the Duistermaat-Heckman theorem. Hence, $\omega()$
defines a class in $\cald((M_0^a,\omega_0^a),[\nu_0^a],\ell)$ which
we define to be $\Theta(a)$.

To see that $\Theta$ is well-defined,
suppose that $\omega'$ is a symplectic form on $M$ equivalent to $\omega$.
Let $F$ be a self-diffeomorphism of $M$ realizing the equivalence, so
$F^*\omega()=\omega'()$. The map
$F$ descends to a self-diffeomorphism $\bar F$ of $[0,\ell]\times M_0^a$ commuting with the
projection onto $[0,\ell]$. Hence, $\bar F$ is of the form
$\bar F(\lambda,x)=(\lambda,\bar F_\lambda(x))$ where $\bar F_\lambda$ is a self-diffeomorphism
of $M_0^a$ such that $\bar F_\lambda^*\omega(\lambda)=\omega'(\lambda)$
and $\bar F_0={\rm id}$.
For $t\in[0,1]$, let $\omega_t\:[0,\ell]\to\Omega^{sym}(M_0^a)$ be defined by
$\omega_t(\lambda)=\bar F_{t\lambda}^*\omega$. The map $t\mapsto \omega_t()$ is
a path in $\Omega^{sym}(M_0^a)$ from $\omega()$ to $\omega'()$. This shows that
the two forms are cohomologuous and so
$\Theta$ is well-defined.

Let us now prove that $\Theta$ is surjective. Let $\omega()$ represent a
class in $\cald((M_0^a,\omega_0^a),[\nu_0^a],\ell)\big)$.
Let $\bbs_0^a\to M_0^a$ be the $S^1$-bundle associated to $\nu_0^a$.
Using the normal form for reduced spaces \cite[\S\,30.3]{CS}, we  can
extend the map $\omega()$
to a smooth map $\omega\:[-\varepsilon,1+\varepsilon]\to\Omega^{sym}(M_0^a)$.
For such map there is a symplectic form $\tilde\omega$ on $\bbs_0^a\times[-\varepsilon,1+\varepsilon]$ such that the $S^1$-action is Hamiltonian
with moment map the projection onto $[-\varepsilon,1+\varepsilon]$, as shown in \cite[Proposition~5.8]{MDS}.
Performing symplectic cuts at $0$ and $1$ provides a \simple\ manifolds $(N;M_0^a,N_1)$
defining a class $a\in \cals^0([\nu_0^a],\omega_0^a,1,\ell)$  and using $\alpha={\rm id}$
so that $\Theta(a)=[\omega()]$.

To prove the injectivity of $\Theta$, suppose $a,a'\in\cals^0([\nu_0^a],\omega_0^a,1,\ell)$
are represented by symplectic forms $\tilde\omega$ and $\tilde\omega'$ on $M=\calc_\chi(M_0,\nu_0,\ell)$.
These give rise to
$\omega()$ and $\omega'()$ in $\cald((M_0^a,\omega_0^a),[\nu_0^a],\ell)$
representing $\Theta(a)$ and $\Theta(a')$.
If $\Theta(a)=\Theta(a')$, there exists a path
$\omega_t()\in\cald((M_0^a,\omega_0^a),[\nu_0^a],\ell)$ joining $\omega()$ to $\omega'()$.
Because of the cohomology constraint in the definition of $\cald((M_0^a,\omega_0^a),[\nu_0^a],\ell)$,
the cohomology class of $\omega_t(\lambda)$ is independent of $t$. By Moser's theorem
\cite[Theorem~7.3]{CS}, there exists an isotopy $\rho_t\:M_0^a\times[0,\ell]\to M_0^a\times[0,\ell]$, with $\rho_0={\rm id}$,
such that $\omega_t()=\rho_t^*\omega()$. This isotopy may be covered
by an isotopy $\tilde\rho_t\:M\to M$ with $\tilde\rho_0={\rm id}$.
Let $\tilde\omega_t=\tilde\rho_t^*\tilde\omega$.
By \cite[Proposition~5.8]{MDS}, we may deduce that $\tilde\omega_1=\tilde\omega'$.
This proves that $a=a'$, completing the proof.
\cqfd

\thref{T.symcod1} reduces the identification of $\cals^0([\nu_0^a],\omega_0^a,1,\ell)$
to computing $\pi_0\big(\cald((M_0^a,\omega_0^a),[\nu_0^a],\ell)\big)$.
We only have results when the latter is reduced to one element.

\begin{Theorem}\label{T.S11conn}
Suppose that $r_1=1$. Then
$\pi_0\big( \cald((M_0^a,\omega_0^a),[\nu_0^a],\ell)\big) = *$ if 
$[\omega_0^a]$ and $e(\nu_0^a)$ are linearly dependent in the vector space $H^2_{dr}(M_0^a)$ of  de Rham cohomology.
\end{Theorem}

The linear dependence condition is automatically fulfilled when $H^2_{dr}(M_0^a)\approx\bbr$,
as when  $M_0^a$ is a complex Grassmannian or $\tilde G_2(\bbr^{m+2})$ of Example~1.1.6.

\sk{1}

\preu
Let $\omega\:[0,\ell]\to\Omega^{sym}(M_0^a)$ represent an element
of $\cald((M_0^a,\omega_0^a),[\nu_0^a],\ell)$. As $[\omega_0^a]\neq 0$, 
our hypothesis of linear dependence implies that there is 
a unique $s\in\bbr$ such that $e(\nu_0^a)=s\,[\omega_0^a]$. 
Hence, 
$$
[\omega(\lambda)] = [\omega_0^a] + \lambda\,e(\nu_0^a) = (1+\lambda s)[\omega_0^a] \, .
$$
As 
$[\omega(\lambda)]\neq 0$, we know that $(1+\lambda s)>0$. The symplectic
form $(1+\lambda s)^\mun\omega(\lambda)$ thus satisfies
$[(1+\lambda s)^\mun\omega(\lambda)]=[\omega_0^a]$. By Moser's theorem
\cite[Theorem~7.3]{CS}, there exists an isotopy
$$\rho_\lambda\:M_0^a\to M_0^a,$$ with $\rho_0={\rm id},$
such that $\omega(\lambda)=\rho_\lambda^*\omega_0^a$.
Hence, the formula
$$
\omega_t(\lambda)=(1+\lambda s)\rho_{t\lambda}^*\omega_0^a \quad (t\in[0,1])
$$
defines a path in $\cald((M_0^a,\omega_0^a),[\nu_0^a],\ell)$ joining
$\omega$ to $(1+\lambda s)\omega_0^a$. This shows that
$\pi_0\big( \cald((M_0^a,\omega_0^a),[\nu_0^a],\ell)\big)$ has
only one element.
\cqfd

Using \thref{T.symcod1}, \thref{T.S11conn} and its proof have the following corollary.

\begin{Corollary}\label{T.S1=1}
Let $(M,M_0,M_1)$ be a $[\nu_0^a]$-simple Hamiltonian $T$-manifold with $T$-size $\ell$ and associated character $\chi$.
Suppose that $r_0=r_1=1$ and that $e(\nu_0^a)=s\, [\omega_0^a]$
for some $s\in\bbr$.
Then there exits a $T$-equivariant symplectomorphism $\alpha\:\calc_\chi(M_0^a,\nu_0^a,\ell)\fl{\approx}M$
such that $\alpha|_{M_0}=\alpha_0$. Moreover,
$(M_1,\omega_1)$ is symplectomorphic to $(M_0^a,(1+s\ell)\,\omega_0^a)$. \cqfd
\end{Corollary}

As a corollary below, we may reproduce Delzant's result \cite[Theorem 1.2]{De} in a slightly  
more precise way, 
with essentially the same proof rephrased in our framework.
For the diagonal action of $S^1$ on $\bbc^{m+1}$, with moment map $\tilde\Phi(z)=\frac{1}{2}|z|^2$, 
denote by $(\bbc P^m)_\ell$ the symplectic reduction at $\ell$:
$$
(\bbc P^m)_\ell = \bbc^{m+1}\mathop{\,\big/\!\!\big/\,}_{\!\textstyle \ell\kern 7pt}S^1 \, .
$$
We also consider  the symplectic cut $(\widehat{\bbc P^m)}_\ell$ of $\bbc^{m+1}$ at $\ell$,
equipped with the induced $S^1$-action and induced moment map $\hat\phi\:(\widehat{\bbc P^m)}_\ell\to [0,\ell]$.
Observe that $(\widehat{\bbc P^m)}_\ell$ is symplectomorphic to $(\bbc P^{m+1})_\ell$.
Indeed, as the symplectic forms vary linearly in $\ell$, it is enough to prove this for $\ell=1$. 
But $(\widehat{\bbc P^m)}_1$ and $(\bbc P^{m+1})_1$ are both toric manifolds admitting as
moment polytope an $(m+1)$-simplex intersecting the weight lattice at its vertices.

\begin{Corollary}\label{T.Sprojsp}
Let $(M^{2m},M_0,M_1)$ be a \simple\ $S^1$-manifold of $S^1$-size $\ell$,
with $M_0$ a single point. Then 
\begin{enumerate}
\item $M$ is $S^1$-equivariantly symplectomorphic to $(\bbc P^m)_\ell$, endowed with a standard $S^1$-action
(multiplication on a single coordinate).
\item $M_1$ is symplectomorphic to $(\bbc P^{m-1})_\ell$.
\end{enumerate}
\end{Corollary}

\preu Let $(W,W_0,W_1)=(\widehat{\bbc P^m)}_\ell,pt,(\bbc P^m)_\ell)$ and, 
for $A\subset\bbr$, let $X_A=\{z\in\bbc\mid |z|\in A\}$. Let $0<\varepsilon<\varepsilon'<\ell$. 
Performing a symplectic cut to $W$ at $\varepsilon$ gives rise to 
two simple manifolds, the ``lower'' one  $(W_-,pt,V_\epsilon)$ and the ``upper'' one $(W_+,V_\epsilon,W_1)$, 
together with symplectic $S^1$- equivariant embeddings $h_-\:X_{(0,\varepsilon)}\to W_-$  and 
$h_+\:X_{\varepsilon,\varepsilon'}\to W_+$. From these, one can recover $W$.
The quotient map $p\:X_{[0,\ell]}$ induces to an $S^1$- equivariant symplectomorphism
\begin{equation}\label{T.Sprojsp-eq50}
W\approx \big(W_--V_\varepsilon\big)  \cup_{h_-} X_{(0,\varepsilon')} \cup_{h_+}
\big(W_+-V_\varepsilon\big) \, .
\end{equation}
In the same way, performing a symplectic cut of $M$ at $\varepsilon$ gives rise to 
two simple manifolds $(M_-,pt,N_\epsilon)$ and $(M_+,N_\epsilon,M_1)$. If $\Phi\:M\to[0,\ell]$ denotes
the moment map, we get $S^1$-symplecto\-morphisms $\alpha_-\:\Phi^\mun([0,\varepsilon))\to M_--N_\varepsilon$
and $\alpha_+\:\Phi^\mun((\varepsilon,\ell])\to M+-N_\varepsilon$.
By the local forms around a fixed point, there is an $S^1$- equivariant symplectomorphism $q\:U\to U'$ between  
neighborhoods $U$ and $U'$ of $W_0$ and $M_0$ respectively. Choose $\varepsilon'$ small enough
so that $p(X_{(0,\varepsilon')})\subset U$. We thus get two symplectic $S^1$- equivariant embeddings
$$
g_-=\alpha_-\pcirc q\pcirc p\:X_{(0,\varepsilon)}\to W_-\ \  \mbox{ and }\ \ 
g_+=\alpha_+\pcirc q\pcirc p\:X_{\varepsilon,\varepsilon'}\to W_+,
$$
and an $S^1$- equivariant symplectomorphism
\begin{equation}\label{T.Sprojsp-eq50b}
M\approx \big(M_--N_\varepsilon\big)  \cup_{g_-} X_{(0,\varepsilon')} \cup_{g_+}
\big(M_+-N_\varepsilon\big) \, .
\end{equation}
The symplectomorphism $q$ induces an $S^1$- equivariant symplectomorphism $q_-\:W_-\to M_-$ such that
$q_-\pcirc h_-=g_-$, and also a symplectomorphism $q_\varepsilon\:V_\varepsilon\to N_\varepsilon$. 

In order to get an $S^1$- equivariant symplectomorphism from $W$ to $M$
it is then enough, given~\eqref{T.Sprojsp-eq50} and~\eqref{T.Sprojsp-eq50b},
to construct an $S^1$- equivariant symplectomorphism $q_+\:W_+\to M_+$ such that
$q_+\pcirc h_+=g_+$. 

The problem may be reformulated as follows. 
Let $Y$ be the upper manifold of the symplectic cut of $X_{(0,\varepsilon')}$ at $\varepsilon$. 
The embedding $h_+$ extends to a symplectic $S^1$- equivariant embedding $\hat h_+\:Y\to W_+$
onto a tubular neighbourhood of $V_\varepsilon$ in $W_+$; $h_+$ and $\bar h_+$ determine each other.
In the same way, $g_+$ extends to  a symplectic $S^1$- equivariant embedding $\hat g_+\:Y\to M_+$
onto a tubular neighbourhood of $N_\varepsilon$ in $M_+$; $g_+$ and $\bar g_+$ determine each other.
We are then looking for an $S^1$- equivariant symplectomorphism $q_+\:W_+\to M_+$ such that
$q_+\pcirc \hat h_+=\hat g_+$.  As $\hat g_+$ coincides with $q_\varepsilon$ on $V_\varepsilon$,
it actually suffices to construct an $S^1$- equivariant symplectomorphism $q_+\:W_+\to M_+$ extending $q_\varepsilon$.
Indeed, by the uniqueness of $S^1$-invariant tubular neighbourhood of $N_\epsilon$ up to symplectomorphism,
it will be possible, taking $\varepsilon'$ smaller if necessary, to modify $q_+$ by an isotopy so that    
the condition $q_+\pcirc \hat h_+=\hat g_+$ remains true.

By construction, $W_+$ is identified with 
$(\calc_{\rm id}((\bbc P^{m-1})_\varepsilon,\eta,\ell-\varepsilon))$, where $\eta$ is the Hopf bundle.
The simple manifold $(M_+,N_\epsilon,M_1)$ has $S^1$-size $\ell-\varepsilon$ and the existence of
the diffeomorphism $\hat g_+$ implies that $g_\varepsilon^*(\nu(N_\varepsilon))=\eta$.
By \corref{T.S1=1}, $g_\varepsilon$ extends to an $S^1$-equivariant symplectomorphism $q_+\:W_+\to M_+$ as required.
\cqfd

\section{Examples of polygon spaces}\label{se:poly}

This section  provides examples using polygon spaces. We recall below some minimal theory to
state the results. For more developments, classification and references, 
see e.g. \cite{HK} and \cite{Ha}.  

Let $\alpha =(\alpha _1,\dots ,\alpha _{n}) \in \rmcr$, where
$\rmcr:=\{(\alpha_1,\dots,\alpha_n)\in\bbr^n\mid 0< \alpha_1\leq\cdots\leq \alpha_n\}$.
Let $S^2_{\alpha_i}$ denote the sphere in $\bbr^3$ with radius $\alpha_i$. We
identify $\bbr^3$ with $so(3)^*$ so that the Lie-Kirillov-Kostant-Souriau
symplectic structure gives $S^2_{\alpha_i}$ the symplectic volume $2\alpha_i$.

\begin{definition}\label{def:polygon}
The \textbf{polygon space} $\pol{\alpha}$ is the symplectic reduction at $0$
 $$ 
\pol{\alpha} = \bigg(\prod_{i=1}^m S^2_{\alpha _i}\bigg)
 \mathop{\,\bigg/\!\!\!\bigg/\,}_{\!\textstyle 0\kern 7pt}SO_3
$$
for the the diagonal co-adjoint action of $SO(3)$. 
\end{definition} 

The moment map for the co-adjoint action on the product of spheres
maps 
$\rho\mapsto\sum\rho_i$, so
we get
\begin{equation}\label{E-conf}
 \pol{\alpha}=\bigg\{\rho=(\rho _1,\dots,\rho _m)\in (\bbr^3)^m \biggm |
 \forall i, |\rho _i|=\alpha_i
 \hbox{ and } \sum_{i=1}^m\rho _i=0  \bigg\}\biggm/SO_3
\end{equation}
 as the moduli space of spatial configurations of a polygon with length-side vector $\alpha$.
Note that $\pol{\alpha}$ is denoted by ${\rm Pol\,}(\alpha)$ in \cite{HK} and by $\nua{n}{3}(\alpha)$ in \cite{Ha}). 
The origin is a regular value for the moment map if and only if there is no aligned configuration, that
is the equation 
$$
\sum_{i=1}^n\epsilon_i\alpha_i = 0 
$$
has no solution with $\epsilon_i=\pm 1$. Such length vectors $\alpha$ are called \textbf{generic}.

When $\alpha_i\neq\alpha_j$ for some $i,j$, then  $\Phi_{i,j}(\rho)=|\rho_i+\rho_j|$
defines a smooth function $\Phi_{i,j}\:\pol{\alpha}\to \bbr$. This is the moment map
of a Hamiltonian $S^1$-action on $\pol{\alpha}$, a particular case of a \textbf{bending flow} \cite{KM}.
It acts on $\rho$ by rotating $\rho_i$ and $\rho_j$ at constant speed around the axis $\rho_i+\rho_j$.
The critical points for $\Phi_{i,j}$ are those configurations $\rho$ for which $\{\rho_k\mid k\neq i,j\}$
generate a one-dimensional space.

If $\alpha\in\rmcr$ satisfies the inequalities
\begin{equation}\label{chanbre1}
\alpha_n< \sum_{i<n} \alpha_i   \ \hbox{ and } \alpha_n+ \alpha_1 > \sum_{i=2}^{n-1} \alpha_i \, ,
\end{equation}
then $\pol{\alpha}$ is known to be diffeomorphic to $\bbc P^{n-3}$, as shown in \cite[Example~2.6]{Ha}. 
Using \corref{T.Sprojsp}, we get a precise symplectic description.

\begin{Proposition}\label{P.chambre1}
Let $\alpha =(\alpha _1,\dots ,\alpha _{n}) \in \rmcr$ ($n\geq 4$) satisfying~\eqref{chanbre1}. 
Then $\pol{\alpha}$ is symplectomorphic to $(\bbc P^{n-3})_\ell$ for 
$\ell=\alpha_1+\cdots+\alpha_{n-1}-\alpha_n$.
\end{Proposition}

\preu
Since $n\geq 4$, the second equation in~\eqref{chanbre1} implies that
\begin{equation}\label{chanbre2}
\alpha_n-\alpha_{n-1} > \alpha_2 + \cdots \alpha_{n-2} - \alpha_1  > 0 \, .
 \end{equation}
Hence, the bending flow $\Phi=\Phi_{n,n-1}$ is defined, with image 
$$
I=[\alpha_n-\alpha_{n-1},\alpha_1 + \cdots \alpha_{n-2}] \, ,
$$ 
an interval of length $\ell=\alpha_1+\cdots+\alpha_{n-1}-\alpha_n$. The fact that $\alpha\in \rmcr$
together with the second inequality of~\eqref{chanbre2} imply that there are no critical points
for $\Phi$ in the interior of $I$. Hence, $\Phi$ makes $\pol{\alpha}$ a simple
Hamiltonian manifold with $S^1$-size equal to $\ell$. The manifold $\Phi^\mun(\alpha_1 + \cdots \alpha_{n-2})$ 
is equal to a point. \proref{P.chambre1} then follows from \corref{T.Sprojsp} (exchanging the role of
$M_0$ and $M_1$).
\cqfd

We now study the operation of adding a tiny edge to a polygon.
Let $\alpha =(\alpha _1,\dots ,\alpha _{n}) \in \rmcr$ 	be generic.
If $\varepsilon>0$ is small enough, then, for all integer $j\in\{1,\dots,n\}$, the
$n$-tuple
$$
\alpha(j,\delta)= (\alpha_1,\dots,\alpha_{j-1},\alpha_j+\delta,\alpha_{j+1},\dots,\alpha_n)
$$
belongs to $\rmcr$ and is generic when $|\delta|\leq\varepsilon$. We say that 
$\varepsilon$ is \textbf{$\alpha$-tiny}. The manifolds $\pol{\alpha(j,\delta)}$ are then 
canonically diffeomorphic to $\pol{\alpha}$, see \cite[Lemma~1.2 and its proof]{Ha}.

We shall now describe the symplectic manifold $\pol{\alpha^\varepsilon}$ where
\begin{equation}\label{E-alphaepsi}
 \alpha^\varepsilon=(\varepsilon,\alpha _1,\dots ,\alpha _{n})\in\bbr_{\scriptscriptstyle\nearrow}^{n+1}
\end{equation}
and $\varepsilon$ is $\alpha$-tiny. For convenience, we will now index the coordinates by $0$ to $n$.
We check that the bending flow 
$$
\Phi_ {j,0}\: \pol{\alpha^\varepsilon} \to I_j= [\alpha_j-\varepsilon,\alpha_j+\varepsilon]
$$ 
is well-defined and makes $\pol{\alpha^\varepsilon}$ a simple Hamiltonian $S^1$-manifold of 
$S^1$-size equal to $2\varepsilon$, with 
$M_0=\pol{\alpha(j,-\varepsilon)}$ and $M_1=\pol{\alpha(j,\varepsilon)}$.

For $i=0,\dots n$, consider the space $E_i$ of configurations $\rho$ as in~\eqref{E-conf} such that  
$\rho_i=(0,0,\alpha_i)$. This is the total space of a principal $S^1$-bundle $\xi$ over 
$\pol{\alpha^\varepsilon}$, or over $\pol{\alpha(i,-\varepsilon)}$ if $1\leq i \leq n$.
We also denote by $\xi_i$ its associated complex line bundle.

\begin{Proposition}\label{P.chambre2}
Let $\alpha\in \rmcr$ and let $\varepsilon>0$ be tiny for  $\alpha$.
Then for each $1\leq j\leq n$, the bending flow $\Phi_ {j,0}$ makes 
the manifold $\pol{\alpha^\varepsilon}$ $S^1$-equivariantly symplectomorphic
to $\calc_{\rm id}(\pol{\alpha(j,-\varepsilon)},\xi_j,2\varepsilon)$. 
\end{Proposition}

\noindent For two descriptions of $\pol{\alpha^\varepsilon}$ as a smooth manifold,
see \cite[Proposition~2.2]{Ha}.

\medskip

 \preu
 Choose an $S^1$-invariant almost complex structure on $\pol{\alpha^\varepsilon}$ compatible with
 the symplectic form and let $\nu$ be the normal bundle to $\pol{\alpha(j,-\varepsilon)}$
 in $\pol{\alpha^\varepsilon}$. 
As $r_0=r_1=1$, \corref{T.S1=1} implies that there is an $S^1$-equivariant symplectomorphism
from $\calc_{\rm id}(\pol{\alpha(j,-\varepsilon)},\nu,2\varepsilon)$ to $\pol{\alpha^\varepsilon}$.
We have to identify $\nu$ with $\xi_j$.

The symplectic reduction of $\pol{\alpha^\varepsilon}$ at $\lambda\in I_j$ is
$\pol{\alpha(j,\lambda)}$. Identifying the latter with $\pol{\alpha}$ gives a symplectic form 
$\omega_\lambda\in \Omega^2(\pol{\alpha})$ which, by the Duistermaat-Heckman theorem
satisfies the equation
$$
[\omega(\lambda)]=[\omega_0^a]+\lambda e(\nu)
$$
in $H^2_{dr}(\pol{\alpha})$. Hence,
\begin{equation}\label{P.chambre2-eq10}
e(\nu) = \frac{d}{\lambda}[\omega(\lambda)] \, .
 \end{equation}
But, by \cite[Remark~7.5.d]{HK},
\begin{equation}\label{P.chambre2-eq20}
 [\omega(0)] = \sum_{i=1}^n \alpha_i e(\xi_i) \, .
\end{equation}
By \eqref{P.chambre2-eq10} and \eqref{P.chambre2-eq10}, we deduce that $\nu$ is isomorphic to $\xi_j$.
\cqfd

\begin{Proposition}\label{P.chambre3}
Let $\alpha=(\alpha_0,\dots,\alpha_n)\in \bbr_{\scriptscriptstyle\nearrow}^{n+1}$ satisfying
\begin{equation}\label{E-chanbre3}
\alpha_n+\alpha_0< \sum_{i<n} \alpha_i   \ \hbox{ and } \alpha_n+ \alpha_1 > \sum_{i=2}^{n-1} \alpha_i \, ,
\end{equation}
let $\ell=\alpha_0+\cdots+\alpha_{n-1}-\alpha_n$.
Then $\pol{\alpha}$ is symplectomorphic to a symplectic cut of $(\bbc P^{m-2})_\ell$ so that the
symplectic slice has size $\ell-2\alpha_0$.
\end{Proposition}

In particular, $\pol{\alpha}$ is diffeomorphic to $\bbc P^{m-2}\, \sharp\, \overline{\bbc P}^{m-2}$.
For a generalization of this fact, see \cite[Example~2.12]{Ha}.

\medskip

\preu
We note that $\alpha=\beta^{\alpha_0}$ in the sense of~\eqref{E-alphaepsi}, where $\beta=(\alpha_1,\dots,\alpha_n)$
satisfies~\eqref{chanbre1}. 
We use \proref{P.chambre2} and its notations, with the bending flow
$\Phi_{n,0}$. Hence, $\pol{\alpha}$ is symplectomorphic to  
$\calc_{\rm id}(\pol{\alpha(n,-\alpha_0)},\xi_n,2\alpha_0)$.  Using~\eqref{P.chambre2-eq10}
and \cite[Proposition~7.3]{HK}, we deduce that $e(\xi_n)=-1$.

Thus, $\Phi_{n,0}$ makes $\pol{\alpha}$ a simple Hamiltonian manifold $(\pol{\alpha},M_0,M_1)$ with
$M_0=(\bbc P^{n-3})_{\ell}$ and $M_1=(\overline{\bbc P}^{n-3})_{\ell-2\alpha_0}$, using
\proref{P.chambre1} to identify $M_0$.
\cqfd

\vskip5mm

\noindent

\begin{minipage}{8.5cm}
\footnotesize
Jean-Claude HAUSMANN\\
Math\'ematiques-Universit\'e\\ B.P. 24\\
CH-1211 Gen\`eve 4, Switzerland\\[2pt]
Jean-Claude.Hausmann@unige.ch
\end{minipage}

\vskip 5mm

\begin{minipage}{6.5cm}
\footnotesize
Tara HOLM\\
Department of Mathematics\\
Cornell University\\
Ithaca, NY  14853-4201, USA\\[2pt]
tsh@math.cornell.edu
\end{minipage}

\end{document}